\documentclass[11pt,bezier]{article}
\usepackage{amsmath}
\usepackage{amsfonts,amsthm,amssymb}
\usepackage{amsfonts}
\usepackage{graphics}
\usepackage{amssymb}
\newtheorem{lem}{Lemma}[section]
\newtheorem{thm}[lem]{Theorem}

\newtheorem{conj}{Conjecture}

\theoremstyle{definition}

\begin{document}
\title{Connectivity keeping stars or double-stars in 2-connected graphs
\footnote{The research is supported by NSFC (Nos.11401510, 11531011) and NSFXJ(No.2015KL019).}}
\author{Yingzhi Tian$^{a}$ \footnote{Corresponding author. E-mail: tianyzhxj@163.com (Y.Tian),
mjx@xju.edu.cn (J.Meng), hjlai@math.wvu.edu (H. Lai), 200661000016@jmu.edu.cn (L. Xu).}, Jixiang Meng$^{a}$, Hong-Jian Lai$^{b}$, Liqiong Xu$^{c}$ \\
{\small $^{a}$College of Mathematics and System Sciences, Xinjiang
University, Urumqi, Xinjiang 830046, PR China}\\
{\small $^{b}$Department of Mathematics, West Virginia University,
Morgantown, WV 26506, USA}\\
{\small $^{c}$School of Science, Jimei University, Xiamen, Fujian 361021, PR China}}

\date{}

\maketitle

\noindent{\bf Abstract } In [W. Mader, Connectivity keeping paths in $k$-connected graphs, J. Graph Theory 65 (2010) 61-69.], Mader conjectured that for every positive integer $k$ and every finite tree $T$ with order $m$,  every $k$-connected, finite graph $G$ with $\delta(G)\geq \lfloor\frac{3}{2}k\rfloor+m-1$ contains a subtree $T'$ isomorphic to $T$
such that $G-V(T')$ is $k$-connected. In the same paper, Mader proved that the conjecture is true when $T$ is a path. Diwan and Tholiya [A.A. Diwan, N.P. Tholiya, Non-separating trees in connected graphs, Discrete Math. 309 (2009) 5235-5237.] verified the conjecture when $k=1$.  In this paper, we will prove that Mader's conjecture is true when $T$ is a star or double-star and $k=2$.

\noindent{\bf Keywords:} 2-Connected graphs; Stars; Double-stars; Mader's Conjecture

\section{Introduction}

In this paper, $graph$ always means a finite, undirected graph without multiple edges and without loops. For graph-theoretical terminologies and
notation not defined here, we follow \cite{Bondy}. For a graph $G$, the vertex set, the edge set, the minimum degree and the connectivity number of $G$
are denoted by $V(G)$, $E(G)$, $\delta(G)$ and $\kappa(G)$, respectively. The $order$ of a graph $G$ is the cardinality of its vertex set, denoted by $|G|$. $k$ and $m$ always denote positive integers.

In 1972, Chartrand, Kaugars, and Lick proved the following well-known result.

\begin{thm}{\cite{Chartrand}} Every $k$-connected graph $G$ of minimum degree
$\delta(G)\geq \lfloor\frac{3}{2}k\rfloor$ has a vertex $u$ with $\kappa(G-u)\geq k$.
\end{thm}

Fujita and Kawarabayashi proved in \cite{Fujita} that every $k$-connected graph $G$ with minimum degree at least $\lfloor\frac{3}{2}k\rfloor+2$ has an edge $e$ such that $G-V(e)$ is still $k$-connected. They conjectured that there are similar results for the existence of connected subgraphs of prescribed order $m\geq3$ keeping the connectivity.

\begin{conj}{\cite{Fujita}} For all positive integers $k, m$, there
is a (least) non-negative integer $f_k(m)$ such that every $k$-connected graph $G$ with $\delta(G)\geq \lfloor\frac{3}{2}k\rfloor-1+f_k(m)$ contains a connected subgraph $W$ of exact order $m$ such that $G-V(W)$ is still $k$-connected.
\end{conj}

They also gave examples in \cite{Fujita} showing that $f_k(m)$ must be at least $m$ for all positive integers $k, m$. In \cite{Mader1}, Mader proved that $f_k(m)$ exists and $f_k(m)=m$ holds for all $k, m$.

\begin{thm}{\cite{Mader1}}
Every $k$-connected graph $G$ with $\delta(G)\geq\lfloor\frac{3}{2}k\rfloor+m-1$
for positive integers $k, m$ contains a path $P$ of order $m$ such that $G-V(P)$ remains $k$-connected.
\end{thm}

In the same paper, Mader \cite{Mader1} asked whether the result is true for any other tree $T$ instead of a path, and gave the following conjecture.

\begin{conj}{\cite{Mader1}} For every positive integer $k$ and every finite tree $T$, there is a least non-negative integer $t_k(T)$, such that every $k$-connected, finite graph $G$ with $\delta(G)\geq \lfloor\frac{3}{2}k\rfloor-1+t_k(T)$ contains a subgraph $T'\cong T$ with $\kappa(G-V(T'))\geq k$.
\end{conj}

Mader showed that $t_k(T)$ exists in \cite{Mader2}.

\begin{thm}{\cite{Mader2}}
Let $G$ be a $k$-connected graph with $\delta(G)\geq 2(k-1+m)^2+m-1$ and let
$T$ be a tree of order $m$ for positive integers $k, m$. Then there is a tree $T'\subseteq G$ isomorphic to $T$ such that $G-V(T')$ remains $k$-connected.
\end{thm}

Mader further conjectured that $t_k(T)=|T|$.

\begin{conj}{\cite{Mader1}} For every positive integer $k$ and every tree $T$, $t_k(T)=|T|$ holds.
\end{conj}

Theorem 1.2 showed that Conjecture 3 is true when $T$ is a path.   Diwan and Tholiya \cite{Diwan} proved that the conjecture holds when $k=1$.  In the next section, we will verify that Conjecture 3 is true when $T$ is a star and $k=2$. It is proved in the last section that Conjecture 3 is true when $T$ is a double-star and $k=2$.

A $block$ of a graph $G$ is a maximal connected subgraph of $G$ that has no cut vertex. Note that any block of a connected graph of order at least two is 2-connected or isomorphic to $K_2$.

For a vertex subset $U$ of a graph $G$, $G[U]$ denotes the subgraph induced by $U$ and $G-U$ is the subgraph induced by $V(G)-U$. The $neighborhood$ $N_G(U)$ of $U$ is the set of vertices in $V(G)-U$ which are adjacent to some vertex in $U$. If $U=\{u\}$, we also use $G-u$ and $N_G(u)$ for $G-\{u\}$  and $N_G(\{u\})$, respectively. The $degree$ $d_G(u)$ of $u$ is $|N_G(u)|$. If $H$ is a subgraph of $G$, we often use $H$ for $V(H)$. For example, $N_G(H)$, $H\cap G$ and $H\cap U$ mean $N_G(V(H))$, $V(H)\cap V(G)$ and $V(H)\cap U$, respectively. If there is no confusion, we always delete the subscript, for example, $d(u)$ for $d_G(u)$, $N(u)$ for $N_G(u)$,  $N(U)$ for $N_G(U)$ and so on. A $tree$ is a connected graph without cycles. A $star$ is a tree that has exact one vertex with degree greater than one. A $double$-$star$ is a tree that has exact two vertices with degree greater than one.

\section{Connectivity keeping stars in 2-connected graphs}

\begin{thm}
Let $G$ be a 2-connected graph with minimum degree $\delta(G)\geq m+2$, where $m$ is a positive integer. Then for a star $T$ with order $m$, $G$ contains a star $T'$ isomorphic to $T$ such that $G-V(T')$ is 2-connected.
\end{thm}

\noindent{\bf Proof.} If $m\leq3$, then $T$ is a path, and the Theorem holds by Theorem 1.2. Thus we assume $m\geq4$ in the following.

Since $\delta(G)\geq m+2$, there is a star $T'\subseteq G$ with $T'\cong T$. Assume $V(T')=\{u, v_1, \cdots, v_{m-1}\}$ and $E(T')=\{uv_i| 1\leq i\leq m-1\}$. We say $T'$ is a star rooted at $u$ or with root $u$. Let $G'=G-T'$. Let $B$ be a maximum block in $G'$ and let $l$ be the number of components of $G'-B$. If $l=0$, then $B=G'$ is 2-connected. So we may assume that $l\geq1$. Let $H_1,\cdots, H_l$ be the components of $G'-B$ with $|H_1|\geq\cdots\geq|H_l|$.

Take such a star $T'$ so that

({\bf P1}) $|B|$ is as large as possible,

({\bf P2}) $(|H_1|,\cdots,|H_l|)$ is as large as possible in lexicographic order, subject to (P1).

We will complete the proof by a series of claims.

\noindent{\bf Claim 1.} $|N(H_i)\cap B|\leq1$ and $|N(H_i)\cap V(T')|\geq1$ for each  $i\in\{1,\cdots,l\}$.

Since $B$ is a block of $G'$, we have $|N(H_i)\cap B|\leq1$ for each  $i\in\{1,\cdots,l\}$. By $G$ is 2-connected, $|N(H_i)\cap V(T')|\geq1$ for each  $i\in\{1,\cdots,l\}$.

\noindent{\bf Claim 2.} $l=1$.

Assume $l\geq2$. By Claim 1, there is an edge $th$ between $T'$ and $H_1$, where $t\in T'$ and $h\in H_1$. Choose a vertex $x\in H_l$. Since $\delta(G)\geq m+2$ and $|N(H_l)\cap B|\leq1$ (by Claim 1), we have $|N(x)\setminus (B\cup t)|\geq m+2-1-1=m$. Thus we can choose a star $T''\cong T$ with root $x$ such that $V(T'')\cap (B\cup t)=\emptyset$. But then either there is a larger block than $B$ in $G-T''$, or $G-T''-B$ contains a larger component than $H_1$ ($H_1\cup t$ is contained in a component of $G-T''-B$), which contradicts to (P1) or (P2).

\noindent{\bf Claim 3.} $|N(t)\cap B|\leq1$ and $|N(t)\cap H_1|\geq2$ for any vertex  $t\in V(T')$.

Assume $|N(t)\cap B|\geq2$. Choose a vertex $x\in H_1$. Since $\delta(G)\geq m+2$ and $|N(H_1)\cap B|\leq1$, we have $|N(x)\setminus (B\cup t)|\geq m+2-1-1=m$. Thus we can choose a star $T''\cong T$ with root $x$ such that $V(T'')\cap (B\cup t)=\emptyset$. But $G-T''$ has a block containing $B\cup t$ as a subset, which contradicts to (P1). Thus $|N(t)\cap B|\leq1$ holds. By $d(t)\geq m+2$ and $|N(t)\cap B|\leq1$, we have $|N(t)\cap H_1|=d(t)-|N(t)\cap B|-|N(t)\cap T'|\geq m+2-1-(m-1)=2$.

\noindent{\bf Claim 4.} For any edge $t_1t_2\in E(T')$, $|N(\{t_1,t_2\})\cap B|\leq1$ holds.

By contradiction, assume $|N(\{t_1,t_2\})\cap B|\geq2$. Because $|N(t_1)\cap B|\leq1$ and $|N(t_2)\cap B|\leq1$, we can assume that there are two distinct vertices $b_1, b_2\in B$ such that $t_1b_1, t_2b_2\in E(G)$. Choose a vertex $x\in H_1$. Since $\delta(G)\geq m+2$ and $|N(H_1)\cap B|\leq1$, we have $|N(x)\setminus (B\cup \{t_1,t_2\})|\geq m+2-1-2=m-1$. Thus we can choose a star $T''\cong T$ with root $x$ such that $V(T'')\cap (B\cup \{t_1,t_2\})=\emptyset$. But then $G-T''$ has a block containing $B\cup \{t_1,t_2\}$ as a subset, which contradicts to (P1).

Because $|N(H_1)\cap B|\leq1$ and $G$ is 2-connected, we have $|N(T')\cap B|\geq1$. The following claim further shows that $|N(T')\cap B|=1$.

\noindent{\bf Claim 5.} $|N(T')\cap B|=1$.

By contradiction, assume $|N(T')\cap B|\geq2$. If $N(u)\cap B\neq \O$, say $N(u)\cap B=\{u'\}$, then we have $N(\{v_1,\cdots,v_{m-1}\})\cap B\subseteq\{u'\}$ by Claim 4. That is, $N(T')\cap B=\{u'\}$, a contradiction. Thus $N(u)\cap B=\O$. Assume, without loss of generality, that there are two distinct vertices $w$ and $w'$ in $B$ such that $v_1w, v_2w'\in E(G)$. If $N(v_3)\cap B=\O$ or $|N(v_3)\cap\{v_1,v_2\}|\leq1$, then we can choose a star $T''$ with order $m$ and root $v_3$ such that $V(T'')\cap(B\cup\{u,v_1,v_2\})=\O$. But then $B\cup \{u,v_1,v_2\}$ is contained in a block of $G-T''$, contradicting to (P1). Thus we assume $v_3$ is adjacent to a vertex $y$ in $B$ and is adjacent to both $v_1$ and $v_2$. Without loss of generality, assume $y$ is distinct from $w$. Then we can choose a star $T''$ with order $m$ and root $u$ such that $V(T'')\cap(B\cup\{v_1,v_3\})=\O$. But $B\cup \{v_1,v_3\}$ is contained in a block of $G-T''$, contradicting to (P1). Thus $|N(T')\cap B|=1$.

By Claim 5, $|N(T')\cap B|=1$. Assume $N(T')\cap B=\{w\}$. Since $G$ is 2-connected, we have $|N(H_1)\cap B|\geq1$. By Claim 1, $|N(H_1)\cap B|=1$. Assume $N(H_1)\cap B=\{z\}$. Let $P$ be a shortest path from $z$ to $w$ going through $H_1$ and $T'$. Assume $P:=p_1p_2\cdots p_{q-1}p_q$, where $p_1=z$, $p_q=w$ and $p_i\in H_1\cup T'$ for each $i\in\{2,\cdots,q-1\}$. Since $P$ is a shortest path,  $|N(p_i)\cap P|=2$ for each $2\leq i\leq q-1$. By $N(T')\cap B=\{w\}$ and $N(H_1)\cap B=\{z\}$, $N(p_i)\cap B\subseteq\{w,z\}\subseteq V(P)$ for each $2\leq i\leq q-1$. Thus $|N(p_i)\cap (B\cup P)|=2$ and $|N(p_i)\cap (V(G)\setminus(B\cup P))|\geq m$  for each $2\leq i\leq q-1$.
This implies $G-B\cup P$ is not empty. For any vertex $x$ in $G-B\cup P$, we have $|N(x)\cap P|\leq 3$. For otherwise, we can find a path $P'$ containing $x$ from $z$ to $w$ going through $H_1$ and $T'$ shorter than $P$, a contradiction. By $\delta(G)\geq m+2$, $|N(x)\cap (G-B\cup P)|\geq m+2-3=m-1$. Then we can find a star $T''\cong T$ with root $x$ such that $T''\cap(B\cup P)=\O$. But then $B\cup P$ is contained in a block of $G-T''$, a contradiction. The proof is thus complete. $\Box$

\section{Connectivity keeping double-stars in 2-connected graphs}

\begin{lem}
Let $G$ be a graph and $T$ be a double-star with order $m$. If there is an edge $e=uv\in E(G)$ such that $|N(u)\setminus v|\geq \lfloor \frac{m}{2}\rfloor-1$, $|N(v)\setminus u|\geq m-3$ and $|N(u)\cup N(v)\setminus \{u,v\}|\geq m-2$, then there is a double-star $T'\subseteq G$ isomorphic to $T$.
\end{lem}

\noindent{\bf Proof.} By $T$ is a double-star, $m\geq4$. Assume the double-star $T$ is constructed from an edge $e'=u'v'$ by adding $r$ leaves to $u'$ and $s$ leaves to $v'$, where $1\leq r\leq s$ and $r+s=m-2$. Then $1\leq r\leq \lfloor \frac{m}{2}\rfloor-1$ and $\lceil \frac{m}{2}\rceil-1\leq s\leq m-3$. Since $|N(u)\setminus v|\geq \lfloor \frac{m}{2}\rfloor-1$, $|N(v)\setminus u|\geq m-3$ and $|N(u)\cup N(v)\setminus \{u,v\}|\geq m-2$, we can find a double-star $T'\cong T$ in $G$ with center-edge $e=uv$, where $u$ is adjacent to $r$ leaves and $v$ is adjacent to $s$ leaves. $\Box$

The main idea of the proof of Theorem 3.2 is similar to that of Theorem 2.1, with much more complicated and different details.

\begin{thm}
Let $T$ be a double-star with order $m$ and $G$ be a 2-connected graph with minimum degree $\delta(G)\geq m+2$. Then $G$ contains a double-star $T'$ isomorphic to $T$ such that $G-V(T')$ is 2-connected.
\end{thm}

\noindent{\bf Proof.} By $T$ is a double-star, $m\geq4$. If $m=4$, then $T$ is a path, and the Theorem holds by Theorem 1.2. Thus we assume $m\geq5$ in the following.

Since $\delta(G)\geq m+2$, there is a double-star $T'\subseteq G$ with $T'\cong T$. Assume $V(T')=\{u, v, u_1, \cdots, u_{r}, v_1, \cdots, v_{s}\}$ and $E(T')=\{uv\}\cup \{uu_i| 1\leq i\leq r\}\cup\{vv_j| 1\leq j\leq s\}$, where $1\leq r\leq s$ and $r+s=m-2$. We say $T'$ is a double-star with center-edge $uv$. Let $G'=G-T'$. Let $B$ be a maximum block in $G'$ and let $l$ be the number of components of $G'-B$. If $l=0$, then $B=G'$ is 2-connected. So we may assume that $l\geq1$. Let $H_1,\cdots, H_l$ be the components of $G'-B$ with $|H_1|\geq\cdots\geq|H_l|$.

Take such a double-star $T'$ so that

({\bf P1}) $|B|$ is as large as possible,

({\bf P2}) $(|H_1|,\cdots,|H_l|)$ is as large as possible in lexicographic order, subject to (P1).

We will complete the proof by a series of claims.

\noindent{\bf Claim 1.} $|N(H_i)\cap B|\leq1$ and $|N(H_i)\cap T'|\geq1$ for each  $i\in\{1,\cdots,l\}$.

Since $B$ is a block of $G'$, we have $|N(H_i)\cap B|\leq1$ for each  $i\in\{1,\cdots,l\}$. By $G$ is 2-connected, $|N(H_i)\cap T'|\geq1$ for each  $i\in\{1,\cdots,l\}$.

\noindent{\bf Claim 2.} $|H_i|\geq2$ for each  $i\in\{1,\cdots,l\}$.

This claim holds because $|N(h_i)\cap H_i|= d(h_i)-|N(h_i)\cap T'|-|N(h_i)\cap B|\geq m+2-m-1=1$ for any vertex $h_i\in H_i$, where $1\leq i\leq l$.

\noindent{\bf Claim 3.} $l=1$.

Assume $l\geq2$. By Claim 1, there is an edge $th$ between $T'$ and $H_1$, where $t\in T'$ and $h\in H_1$.
By Claim 2, we can choose an edge $xy\in E(H_l)$. Since $\delta(G)\geq m+2$ and $|N(H_l)\cap B|\leq1$ (by Claim 1), we have $|N(x)\setminus (B\cup \{y,t\})|\geq m+2-1-2=m-1$ and $|N(y)\setminus (B\cup \{x,t\})|\geq m+2-1-2=m-1$. Thus, by Lemma 3.1, we can choose a double-star $T''\cong T$ with center-edge $xy$ such that $V(T'')\cap (B\cup t)=\emptyset$. But then either there is a larger block than $B$ in $G-T''$, or $G-T''-B$ contains a larger component than $H_1$ ($H_1\cup t$ is contained in a component of $G-T''-B$), which contradicts to (P1) or (P2).

\noindent{\bf Claim 4.} $|N(t)\cap B|\leq1$ and $|N(t)\cap H_1|\geq2$ for any vertex  $t\in V(T')$.

Assume $|N(t)\cap B|\geq2$. Choose an edge $xy\in E(H_1)$. Since $\delta(G)\geq m+2$ and $|N(H_1)\cap B|\leq1$, we have $|N(x)\setminus (B\cup \{y,t\})|\geq m+2-1-2=m-1$ and $|N(y)\setminus (B\cup \{x,t\})|\geq m+2-1-2=m-1$. Thus, by Lemma 3.1, we can choose a double-star $T''\cong T$ with center-edge $xy$ such that $V(T'')\cap (B\cup t)=\emptyset$. But then $B\cup t$ is contained in a block of $G-T''$, which contradicts to (P1). Thus $|N(t)\cap B|\leq1$ holds for any vertex  $t\in V(T')$. By $d(t)\geq m+2$ and $|N(t)\cap B|\leq1$, we have $|N(t)\cap H_1|=d(t)-|N(t)\cap B|-|N(t)\cap T'|\geq m+2-1-(m-1)=2$.

\noindent{\bf Claim 5.} For any edge $t_1t_2\in E(T')$, $|N(\{t_1,t_2\})\cap B|\leq1$ holds.

By contradiction, assume $|N(\{t_1,t_2\})\cap B|\geq2$. Because $|N(t_1)\cap B|\leq1$ and $|N(t_2)\cap B|\leq1$, we can assume that there are two distinct vertices $b_1, b_2\in B$ such that $t_1b_1, t_2b_2\in E(G)$. Choose an edge $xy\in E(H_1)$. Since $\delta(G)\geq m+2$ and $|N(H_1)\cap B|\leq1$, we have $|N(x)\setminus (B\cup \{y,t_1,t_2\})|\geq m+2-1-3=m-2$ and $|N(y)\setminus (B\cup \{x,t_1,t_2\})|\geq m+2-1-3=m-2$. Thus, by Lemma 3.1, we can choose a double-star $T''\cong T$ with center-edge $xy$ such that $V(T'')\cap (B\cup \{t_1,t_2\})=\emptyset$. But then $G-T''$ has a block containing $B\cup \{t_1,t_2\}$ as a subset, which contradicts to (P1).

\noindent{\bf Claim 6.} For any 3-path $t_1t_2t_3$ in $T'$, $|N(\{t_1,t_2,t_3\})\cap B|\leq1$ holds.

By contradiction, assume $|N(\{t_1,t_2,t_3\})\cap B|\geq2$. Then we have $|N(t_2)\cap B|=0$. For otherwise, if $|N(t_2)\cap B|=1$, then we have $|N(\{t_1,t_2,t_3\})\cap B|\leq1$ by $|N(\{t_1,t_2\})\cap B|\leq1$ and $|N(\{t_2,t_3\})\cap B|\leq1$, a contradiction. Because $|N(t_1)\cap B|\leq1$ and $|N(t_3)\cap B|\leq1$, we can assume that there are two distinct vertices $b_1, b_3\in B$ such that $t_1b_1, t_3b_3\in E(G)$. Choose any edge $xy\in E(H_1)$. Since $\delta(G)\geq m+2$ and $|N(H_1)\cap B|\leq1$, we have $|N(x)\setminus (B\cup \{y,t_1,t_2,t_3\})|\geq m+2-1-4=m-3>\lfloor \frac{m}{2}\rfloor-1$ (By $m\geq5$) and $|N(y)\setminus (B\cup \{x,t_1,t_2,t_3\})|\geq m+2-1-4=m-3$.

If $|N(x)\setminus (B\cup \{y,t_1,t_2,t_3\})|>m-3$ or $|N(y)\setminus (B\cup \{x,t_1,t_2,t_3\})|>m-3$, then by Lemma 3.1, we can choose a double-star $T''\cong T$ with center-edge $xy$ such that $V(T'')\cap (B\cup \{t_1,t_2,t_3\})=\emptyset$. But then $G-T''$ has a block containing $B\cup \{t_1,t_2,t_3\}$ as a subset, which contradicts to (P1). Thus we assume $|N(x)\setminus (B\cup \{y,t_1,t_2,t_3\})|=m-3$ and $|N(y)\setminus (B\cup \{x,t_1,t_2,t_3\})|=m-3$, which imply $|N(x)\cap B|=1$ and $|N(y)\cap B|=1$.
Since $|N(H_1)\cap B|\leq1$, we can assume $N(x)\cap B=N(y)\cap B=\{z\}$. Without loss of generality, assume $z\neq b_1$.

If $N(x)\setminus y\neq N(y)\setminus x$, then $|N(x)\cup N(y)\setminus (B\cup \{x,y,t_1,t_2,t_3\})|\geq m-2$. So we can choose a double-star $T''\cong T$ with center-edge $xy$ disjoint from $B\cup \{t_1,t_2,t_3\}$. But then $G-T''$ contains a larger block than $B$, a contradiction. Thus $N(x)\setminus y=N(y)\setminus x$. Because we choose the edge $xy$ in $H_1$ arbitrarily, we conclude that $H_1$ is a complete graph and each vertex not in $H_1$ is adjacent to all vertices in $H_1$ if it is adjacent to one vertex in $H_1$.
In particular, every vertex $t$ in $T'$ is adjacent to all vertices in $H_1$ by Claim 4 and the vertex $z$ in $B$ is adjacent to all vertices in $H_1$.

Let $t_4h_4$ be an edge of graph $G$, where $t_4\in V(T')\setminus\{t_1,t_2,t_3\}$ and $h_4\in V(H_1)$. Let $h_1$ be a vertex in $H_1$ distinct from $h_4$. Then $t_1h_1,  h_1z\in E(G)$. Thus we can choose a double-star $T''\cong T$ with center-edge $t_4h_4$ disjoint from $B\cup \{t_1,h_1\}$. But then $B\cup \{t_1,h_1\}$ is contained in a block of $G-T''$, contradicting to (P1).

Because $|N(H_1)\cap B|\leq1$ and $G$ is 2-connected, we have $|N(T')\cap B|\geq1$. The following claim further shows that $|N(T')\cap B|=1$.

\noindent{\bf Claim 7.} $|N(T')\cap B|=1$.

By contradiction, assume $|N(T')\cap B|\geq2$. If $N(u)\cap B\neq \O$, say $N(u)\cap B=\{u'\}$, then we have $N(\{u_1,\cdots,u_r,v\})\cap B\subseteq\{u'\}$ by Claim 5 and $N(\{v_1,\cdots,v_s\})\cap B\subseteq\{u'\}$ by Claim 6. That is, $N(T')\cap B=\{u'\}$, a contradiction. Thus $N(u)\cap B=\O$. Similarly, we have $N(v)\cap B=\O$. Since $|N(\{u_1,\cdots,u_r\})\cap B|\leq1$ and $|N(\{v_1,\cdots,v_s\})\cap B|\leq1$ (By Claim 6), we have $|N(T')\cap B|=2$. Assume, without loss of generality, that there are two distinct vertices $w$ and $w'$ in $B$ such that $u_1w, v_1w'\in E(G)$.

We first show that any vertex $x$ in $\{u_1,\cdots,u_r,v_1,\cdots,v_s\}\setminus\{u_1,v_1\}$ has no neighbors in $B$. By contradiction, assume there is a vertex in $\{u_1,\cdots,u_r,v_1,\cdots,v_s\}\setminus\{u_1,v_1\}$, say $v_j$ for some $j\in\{2,\cdots,s\}$ (the case $u_i$ for some $i\in\{2,\cdots,r\}$ can be proved similarly), such that $N(v_j)\cap B=\{w'\}$. If $v_j$ is adjacent to $u$ (or $u_1$), then for any edge $vv'$ ($v'$ is a neighbor of $v$ in $H_1$), we have $|N(v)\setminus (B\cup\{u,u_1,v_j,v'\})|\geq m+2-4=m-2$ (or $|N(v)\setminus (B\cup\{u_1,v_j,v'\})|\geq m+2-3=m-1$) and $|N(v')\setminus (B\cup\{u,v,u_1,v_j\})|\geq m+2-1-4=m-3$ (or $|N(v')\setminus (B\cup\{v,u_1,v_j\})|\geq m+2-1-3=m-2$). By Lemma 3.1, we can find a double-star $T''\cong T$ with center-edge $vv'$ such that $T''$ is disjoint from $B\cup\{u,u_1,v_j\}$ (or $B\cup\{u_1,v_j\}$). But then $G-T''$ contains a larger block than $B$, a contradiction. Thus neither $u$ nor $u_1$ is adjacent to $v_j$. Choose a neighbor $v_j'$ of $v_j$ in $H_1$. Since $|N(v_j)\setminus (B\cup\{u,v,u_1,v_1,v_j'\})|\geq m+2-1-3=m-2$ and $|N(v_j')\setminus (B\cup\{u,v,u_1,v_1,v_j\})|\geq m+2-1-5=m-4\geq\lfloor\frac{m}{2}\rfloor-1$ (By $m\geq5$),  we can find a double-star $T''\cong T$ with center-edge $v_jv_j'$ such that $T''$ is disjoint from $B\cup\{u,v,u_1,v_1\}$. But then $G-T''$ contains a larger block than $B$, a contradiction. Thus we have $N(\{u_1,\cdots,u_r,v_1,\cdots,v_s\}\setminus\{u_1,v_1\})\cap B=\O$.

Let $v_2v_2'\in E(G)$, where $v_2'$ is a neighbor of $v_2$ in $H_1$. Since $\delta(G)\geq m+2$ and $N(v_2)\cap B=\O$, we have $|N(v_2)\setminus(B\cup\{u,v,u_1,v_1,v_2'\})|\geq m+2-5=m-3$ and $|N(v_2')\setminus(B\cup\{u,v,u_1,v_1,v_2\})|\geq m+2-1-5=m-4\geq\lfloor\frac{m}{2}\rfloor-1$ (By $m\geq5$). If $|N(v_2)\setminus(B\cup\{u,v,u_1,v_1,v_2'\})|\geq m-2$, then, by Lemma 3.1, we can find a double-star $T''\cong T$ with center-edge $v_2v_2'$ such that $T''$ avoids $B\cup\{u,v,u_1,v_1\}$. But then $G-T''$ contains a larger block than $B$, a contradiction. Thus assume $|N(v_2)\setminus(B\cup\{u,v,u_1,v_1,v_2'\})|=m-3$, which implies $v_2$ is adjacent to both $u_1$ and $v_1$.  For the edge $uv$, we can verify that $|N(u)\setminus(B\cup\{v,u_1,v_1,v_2\})|\geq m+2-4=m-2$ and $|N(v)\setminus(B\cup\{u,u_1,v_1,v_2\})|\geq m+2-4=m-2$. By Lemma 3.1, we can find a double-star $T''\cong T$ with center-edge $uv$ such that $T''$ avoids $B\cup\{u_1,v_1,v_2\}$. But then $B\cup\{u_1,v_1,v_2\}$ is contained in a block of $G-T''$, contradicting  to (P1). Thus Claim 7 holds.

By Claim 7, $|N(T')\cap B|=1$. Assume $N(T')\cap B=\{w\}$. Since $G$ is 2-connected, we have $|N(H_1)\cap B|=1$ by Claim 1. Let $N(H_1)\cap B=\{z\}$. Let $P$ be a shortest path from $z$ to $w$ going through $H_1$ and $T''$. Assume $P:=p_1p_2\cdots p_{q-1}p_q$, where $p_1=z$, $p_q=w$ and $p_i\in H_1\cup T'$ for each $i\in\{2,\cdots,q-1\}$. Since $P$ is a shortest path, $N(p_i)\cap P=\{p_{i-1},p_{i+1}\}$ for $2\leq i\leq q-1$. Because $\delta(G)\geq m+2$ and $N(p_i)\cap B\subseteq\{w,z\}\subseteq P$ for each $2\leq i\leq q-1$, we know $p_i$ has at least $m$ neighbors not in $B\cup P$, that is, $G-B\cup P$ is not empty. For any vertex $x$ in $G-B\cup P$, we have $|N(x)\cap P|\leq 3$. For otherwise, we can find a path $P'$ containing $x$ from $z$ to $w$ going through $H_1$ and $T''$ shorter than $P$, a contradiction. By $\delta(G)\geq m+2$, $|N(x)\cap (G-B\cup P)|\geq m+2-3=m-1$. Choose an edge $xy$ in $G-B\cup P$. Since $|N(x)\setminus (B\cup P\cup y)|\geq m+2-4=m-2$ and $|N(y)\setminus (B\cup P\cup x)|\geq m+2-4=m-2$, we can find a double-star $T''\cong T$ with center-edge $xy$ such that $T''\cap(B\cup P)=\O$. But then $B\cup P$ is contained in a block of $G-T''$, a contradiction. The proof is thus complete.  $\Box$

\end{document}